\documentclass[draft,11pt,english]{article}
\usepackage{amsfonts,amsmath,amsthm,amscd,amssymb,latexsym}

    \usepackage[T2A]{fontenc}
    \usepackage[cp1251]{inputenc}
    \usepackage[russian]{babel}

\begin{document}


\title{ON INTEGRAL CONDITIONS\\IN THE MAPPING THEORY}

\author{V. Ryazanov, U. Srebro and E. Yakubov}

\date{}

\theoremstyle{plain}
\newtheorem{theorem}{Theorem}[section]
\newtheorem{lemma}{Lemma}[section]
\newtheorem{proposition}{Proposition}[section]
\newtheorem{corollary}{Corollary}[section]
\newtheorem{definition}{Definition}[section]
\theoremstyle{definition}
\newtheorem{example}{Example}[section]
\newtheorem{remark}{Remark}[section]
\newcommand{\keywords}{\textbf{Key words.  }\medskip}
\newcommand{\subjclass}{\textbf{MSC 2000. }\medskip}
\renewcommand{\abstract}{\textbf{Abstract.  }\medskip}
\numberwithin{equation}{section}

\def\Xint#1{\mathchoice
   {\XXint\displaystyle\textstyle{#1}}%
   {\XXint\textstyle\scriptstyle{#1}}%
   {\XXint\scriptstyle\scriptscriptstyle{#1}}%
   {\XXint\scriptscriptstyle\scriptscriptstyle{#1}}%
   \!\int}
\def\XXint#1#2#3{{\setbox0=\hbox{$#1{#2#3}{\int}$}
     \vcenter{\hbox{$#2#3$}}\kern-.5\wd0}}
\def\dashint{\Xint-}

\maketitle

\begin{abstract} It is established interconnections between various
integral conditions that play an important role in the theory of
space mappings and in the theory of degenerate Beltrami equations in
the plane.
\end{abstract}

\medskip




\maketitle

\section{Introduction} In the theory of space mappings and in the
Beltrami equation theory in the plane, the integral conditions of
the following type \begin{equation}\label{eq1} \int\limits_{0}^{1}\
\frac{dr}{rq^{\lambda}(r)}\ =\ \infty,\  \ \ \ \ \ \ \ \ \  \lambda\
\in \ (0,1]\ ,
\end{equation} are often met with where the function $Q$ is given say in
the unit ball ${\Bbb B}^n=\{ x\in{\Bbb R}^n: |x|<1\}$ and $q(r)$ is
the average of the function $Q(x)$ over the sphere $|x|=r$, see e.g.
\cite{AIM}, \cite{BGR$_2$}, \cite{Ch}, \cite{GMSV},
\cite{Le}--\cite{Per}, \cite{RSY$_1$}--\cite{RSY$_4$}, \cite{Zo$_1$}
and \cite{Zo$_2$}.
\medskip

On the other hand, in the theory of  mappings called quasiconformal
in the mean, conditions of the type
\begin{equation}\label{eq2} \int\limits_{{\Bbb B}^n} \Phi
(Q(x))\ dx\  <\ \infty\end{equation} are standard for various
characteristics $Q$ of these mappings, see e.g. \cite{Ah},
\cite{Bi}, \cite{Gol}, \cite{Kr}--\cite{Ku}, \cite{Per}--\cite{Rya}
and \cite{Str}. Here $dx$ corresponds to the Lebsgue measure in
${\Bbb R}^n$, $n\ge 2$.
\medskip

In this connection, we establish interconnections between integral
conditions on the function $\Phi$ and between (\ref{eq1}) and
(\ref{eq2}). More precisely, we give a series of equivalent
conditions for $\Phi$ under which (\ref{eq2}) implies (\ref{eq1}).
It makes possible to apply many known results formulated under the
condition (\ref{eq1}) to the theory of the mean quasiconformal
mappings in space as well as to the theory of the de\-ge\-ne\-ra\-te
Beltrami equations in the plane.

\medskip

\section{On some equivalent integral conditions} In this section we
establish equivalence of integral conditions, see also Section 3 for
one more related condition.\medskip

Further we use the following notion of the inverse function for
mo\-no\-to\-ne functions. For every non-decreasing function
$\Phi:[0,\infty ]\to [0,\infty ] ,$ the {\bf inverse function}
$\Phi^{-1}:[0,\infty ]\to [0,\infty ]$ can be well defined by
setting
\begin{equation}\label{eq5.5CC} \Phi^{-1}(\tau)\ =\
\inf\limits_{\Phi(t)\ge \tau}\ t\ .
\end{equation} As usual, here $\inf$ is equal to $\infty$ if the set of
$t\in[0,\infty ]$ such that $\Phi(t)\ge \tau$ is empty. Note that
the function $\Phi^{-1}$ is non-decreasing, too.

\begin{remark}\label{rmk3.333} Immediately by the
definition  it is evident  that
\begin{equation}\label{eq5.5CCC} \Phi^{-1}(\Phi(t))\ \le\ t\ \ \ \ \
\ \ \ \forall\ t\in[ 0,\infty ]
\end{equation} with the equality in (\ref{eq5.5CCC}) except
intervals of constancy of the function $\Phi(t)$.
\end{remark}

Similarly, for every non-increasing function $\varphi:[0,\infty ]\to
[0,\infty ] ,$ set
\begin{equation}\label{eq5.5C} \varphi^{-1}(\tau)\ =\
\inf\limits_{\varphi(t)\le \tau}\ t\ . \end{equation} Again, here
$\inf$ is equal to $\infty$ if the set of $t\in[0,\infty ]$ such
that $\varphi(t)\le \tau$ is empty. Note that the function
$\varphi^{-1}$ is also non-increasing.

\begin{lemma} \label{pr5.AAA} Let $\psi:[0,\infty ]\to [0,\infty ]$ be a
sense--reversing ho\-meo\-mor\-phism and $\varphi:[0,\infty ]\to
[0,\infty ]$ a monotone function. Then
\begin{equation}\label{eq5.5Cc} [\psi\circ\varphi ]^{-1}(\tau)\ =\
\varphi^{-1}\circ\psi^{-1}(\tau)\ \ \ \ \ \ \ \forall \tau\in
[0,\infty ] \end{equation} and
\begin{equation}\label{eq5.CCC} [\varphi\circ\psi ]^{-1}(\tau)\ \le\
\psi^{-1}\circ\varphi^{-1}(\tau)\ \ \ \ \ \ \ \forall \tau\in
[0,\infty ]
\end{equation} with the equality in (\ref{eq5.CCC}) except a countable
collection of $\tau\in [0,\infty].$ \end{lemma}

\begin{remark}\label{rmk4.A} If $\psi$ is a sense--preserving homeomorphism,
then (\ref{eq5.5Cc}) and (\ref{eq5.5CC}) are obvious for every
monotone function $\varphi .$ Similar notations and statements also
hold for other segments $[a,b]$, where $a$ and $b\in [-\infty
,+\infty ]$, instead of the segment $[0,\infty]$.\end{remark}

{\it Proof of Lemma \ref{pr5.AAA}.} We first prove (\ref{eq5.5Cc}).
If $\varphi$ is non-increasing, then
$$\left[\psi\circ\varphi\right]^{-1}\,(\tau)=\,\inf\limits_{\psi\left(\varphi(t)\right)\geq \tau}t
\,=\,\inf\limits_{\varphi(t)\leq\psi^{-1}(\tau)}t\,=\,\varphi^{-1}\circ\psi^{-1}(\tau)\,.$$
Similarly, if $\varphi$ is non-decreasing, then
$$\left[\psi\circ\varphi\right]^{-1}\,(\tau)=\,
\inf\limits_{\psi\left(\varphi(t)\right)\leq \tau}t
\,=\,\inf\limits_{\varphi(t)\geq\psi^{-1}(\tau)}t\,=\,\varphi^{-1}\circ\psi^{-1}(\tau)\,.$$

Now, let us prove (\ref{eq5.CCC}) and (\ref{eq5.5CC}). If $\varphi$
is non-increasing, then applying the substitution $\eta=\psi(t)$ we
have
$$\left[\varphi\circ\psi\right]^{-1}\,(\tau)=\,\inf\limits_
{\varphi\left(\psi(t)\right)\geq \tau}t\,=\,
\inf\limits_{\varphi(\eta)\geq \tau}\psi^{-1}(\eta)\,=\,
\psi^{-1}\left(\sup\limits_{\varphi(\eta)\geq
\tau}\eta\right)\,\leq\,$$$$\leq\,\psi^{-1}\left(\inf\limits_{\varphi(\eta)\leq
\tau}\eta\right)\,=\,\psi^{-1}\circ\varphi^{-1}\left(\tau\right)\,,$$
i.e., (\ref{eq5.CCC}) holds for all $\tau\in[0,\infty].$ It is
evident that here the strict inequality is possible only for a
countable collection of $\tau\in[0,\infty]$ because an interval of
constancy of $\varphi$ corresponds to every such $\tau .$ Hence
(\ref{eq5.5CC}) holds for all $\tau\in[0,\infty]$ if and only if
$\varphi$ is decreasing.

Similarly, if $\varphi$ is non-decreasing, then
$$\left[\varphi\circ\psi\right]^{-1}\,(\tau)=\,\inf\limits_
{\varphi\left(\psi(t)\right)\leq \tau}t\,=\,
\inf\limits_{\varphi(\eta)\leq \tau}\psi^{-1}(\eta)\,=\,
\psi^{-1}\left(\sup\limits_{\varphi(\eta)\leq
\tau}\eta\right)\,\leq$$$$\leq\,\psi^{-1}\left(\inf\limits_{\varphi(\eta)\geq
\tau}\eta\right)\,=\,\psi^{-1}\circ\varphi^{-1}\left(\tau\right)\
,$$ i.e., (\ref{eq5.CCC}) holds for all $\tau\in[0,\infty]$ and
again the strict inequality is possible only for a countable
collection of $\tau\in[0,\infty].$ In the case, (\ref{eq5.5CC})
holds for all $\tau\in[0,\infty]$ if and only if $\varphi$ is
increasing.\bigskip

\begin{corollary} \label{cor3.C} In particular, if $\varphi:[0,\infty ]\to [0,\infty
] $ is a monotone function and $\psi = j$ where $j(t)=1/t,$ then
$j^{-1}=j$ and \begin{equation}\label{eq5.5Cd} [j\circ \varphi
]^{-1}(\tau)\ =\ \varphi^{-1}\circ j(\tau)\ \ \ \ \ \ \ \forall
\tau\in [0,\infty ]
\end{equation} i.e., \begin{equation}\label{eq5.5Ce} \varphi^{-1}(\tau)\ =\
\Phi^{-1}(1/{\tau})\ \ \ \ \ \ \ \forall \tau\in [0,\infty ]
\end{equation} where $\Phi=1/{\varphi},$ \begin{equation}\label{eq5.C}
[\varphi\circ j ]^{-1}(\tau)\ \le\ j\circ\varphi^{-1}(\tau) \ \ \ \
\ \ \ \forall \tau\in [0,\infty ]  \end{equation} i.e.,  the inverse
function of $\varphi(1/t)$ is dominated by $1/\varphi^{-1},$ and
except a countable collection of $\tau\in [0,\infty ] $
\begin{equation}\label{eq5.5Cb} [\varphi\circ j ]^{-1}(\tau)\ =\
j\circ\varphi^{-1}(\tau)\ .\end{equation} $1/\varphi^{-1}$ is the
inverse function of $\varphi(1/t)$ if and only if the function
$\varphi $ is strictly monotone.
\end{corollary}

Further, the integral in (\ref{eq333F}) is un\-der\-stood as the
Lebesgue--Stieltjes integral and the integrals in (\ref{eq333Y}) and
(\ref{eq333B})--(\ref{eq333A}) as the ordinary Lebesgue integrals.
In (\ref{eq333Y}) and (\ref{eq333F}) we complete the definition of
integrals by $\infty$ if $\Phi(t)=\infty ,$ correspondingly,
$H(t)=\infty ,$ for all $t\ge T\in[0,\infty) .$

\begin{theorem} \label{pr4.1aB} Let $\Phi:[0,\infty ]\to [0,\infty ]$ be a
non-decreasing function. Set \begin{equation}\label{eq333E} H(t)\ =\
\log \Phi(t)\ .\end{equation}

Then the equality \begin{equation}\label{eq333Y}
\int\limits_{\Delta}^{\infty} H'(t)\ \frac{dt}{t}\ =\ \infty
\end{equation} implies the equality \begin{equation}\label{eq333F}
\int\limits_{\Delta}^{\infty} \frac{dH(t)}{t}\ =\ \infty
\end{equation} and (\ref{eq333F}) is equivalent to
\begin{equation}\label{eq333B} \int\limits_{\Delta}^{\infty}H(t)\
\frac{dt}{t^2}\ =\ \infty
\end{equation} for some $\Delta>0,$ and (\ref{eq333B}) is equivalent to
every of the equalities: \begin{equation}\label{eq333C}
\int\limits_{0}^{\delta}H\left(\frac{1}{t}\right)\ {dt}\ =\ \infty
\end{equation} for some $\delta>0,$ \begin{equation}\label{eq333D}
\int\limits_{\Delta_*}^{\infty} \frac{d\eta}{H^{-1}(\eta)}\ =\
\infty
\end{equation} for some $\Delta_*>H(+0),$ \begin{equation}\label{eq333A}
\int\limits_{\delta_*}^{\infty}\ \frac{d\tau}{\tau \Phi^{-1}(\tau
)}\ =\ \infty \end{equation} for some $\delta_*>\Phi(+0).$
\medskip

Moreover, (\ref{eq333Y}) is equivalent  to (\ref{eq333F}) and hence
(\ref{eq333Y})--(\ref{eq333A})
 are equivalent each to other  if $\Phi$ is in addition absolutely continuous.
In particular, all the conditions (\ref{eq333Y})--(\ref{eq333A}) are
equivalent if $\Phi$ is convex and non--decreasing.
\end{theorem}

It is necessary to give one more explanation. From the right hand
sides in the conditions (\ref{eq333Y})--(\ref{eq333A}) we have in
mind $+\infty$. If $\Phi(t)=0$ for $t\in[0,t_*]$, then
$H(t)=-\infty$ for $t\in[0,t_*]$ and we complete the definition
$H'(t)=0$ for $t\in[0,t_*]$. Note, the conditions (\ref{eq333F}) and
(\ref{eq333B}) exclude that $t_*$ belongs to the interval of
integrability because in the contrary case the left hand sides in
(\ref{eq333F}) and (\ref{eq333B}) are either equal to $-\infty$ or
indeterminate. Hence we may assume in (\ref{eq333Y})--(\ref{eq333C})
that $\Delta>t_0$ where $t_0\colon =\sup\limits_{\Phi(t)=0}t$,
$t_0=0$ if $\Phi(0)>0$, and $\delta<1/t_0$, correspondingly.\medskip

{\it Proof.} The equality (\ref{eq333Y}) implies (\ref{eq333F})
because except the men\-tio\-ned special case
$$
\int\limits_{\Delta}^{T}\ d\Psi(t)\ \geq\ \int\limits_{\Delta}^{T}
\Psi^{\prime}(t)\ dt\ \ \ \ \ \ \ \ \forall\ T\in(\Delta,\ \infty)
$$
where
$$
\Psi (t)\ \colon =\ \int\limits_{\Delta}^t\frac{d\,H(\tau)}{\tau}\
,\ \ \ \Psi'(t)\ =\ \frac{H'(t)}{t}\ ,
$$
see e.g. Theorem $IV.7.4$ in \cite{Sa}, p. 119, and hence
$$
\int\limits_{\Delta}^{T}\frac{d\,H(t)}{t}\ \geq\
\int\limits_{\Delta}^{T} H^{\prime}(t)\ \frac{dt}{t}\ \ \ \ \ \ \ \
\forall\ T\in(\Delta,\ \infty)
$$

The equality (\ref{eq333F}) is equivalent to (\ref{eq333B}) by
integration by parts, see e.g. Theorem III.14.1  in \cite{Sa}, p.
102. Indeed, again except the mentioned special case, through
integration by parts we have
$$
\int\limits_{\Delta}^{T}\frac{d\,H(t)}{t}\ -\
\int\limits_{\Delta}^{T} H(t)\ \frac{dt}{t^2}\ =\ \frac{H(T+0)}{T}\
-\frac{H(\Delta -0)}{\Delta}\ \ \ \ \ \ \ \ \forall\ T\in(\Delta,\
\infty)
$$
and, if
$$
\liminf\limits_{t\to\infty}\ \frac{H(t)}{t}\ < \infty\ ,
$$
then the equivalence of (\ref{eq333F}) and (\ref{eq333B}) is
obvious. If
$$
\lim\limits_{t\to\infty}\ \frac{H(t)}{t}\ = \infty\ ,
$$
then (\ref{eq333B}) obviously holds, $\frac{H(t)}{t}\ge 1$ for
$t>t_0$ and
$$
\int\limits_{t_0}^{T}\frac{d\,H(t)}{t} =
\int\limits_{t_0}^{T}\frac{H(t)}{t}\ \frac{d\,H(t)}{H(t)} \ge \log\
\frac{H(T)}{H(t_0)} = \log\ \frac{H(T)}{T} + \log\ \frac{T}{H(t_0)}
\to \infty
$$
as $T\to\infty ,$ i.e. (\ref{eq333F}) holds, too.
\medskip

Now, (\ref{eq333B}) is equivalent to (\ref{eq333C}) by the change of
variables $t\rightarrow 1/t.$\medskip

Next, (\ref{eq333C}) is equivalent to (\ref{eq333D}) because by the
geometric sense of integrals as areas under graphs of the
corresponding integrands
$$
\int\limits_0^{\delta} \Psi(t)\ dt\ =\
\int\limits_{\Psi(\delta)}^{\infty}\Psi^{-1}(\eta)\ d\eta\ +\
\delta\cdot\Psi(\delta)
$$
where $\Psi(t)=H(1/t)$, and because by Corollary \ref{cor3.C} the
inverse function for $H\left(1/t\right)$ coincides with $1/H^{-1}$
at all points except a countable col\-lec\-tion.\medskip

Further, set $\psi(\xi)=\log{\xi}.$ Then $H=\psi\circ\Phi$ and by
Lemma \ref{pr5.AAA} and Remark \ref{rmk4.A}
$H^{-1}=\Phi^{-1}\circ\psi^{-1},$ i.e.,
$H^{-1}(\eta)=\Phi^{-1}(e^{\eta}),$ and  by the substitutions
$\tau=e^{\eta},$\,\,$\eta\,=\,\log\ \tau$ we have the equivalence of
(\ref{eq333D}) and (\ref{eq333A}).\medskip

Finally, (\ref{eq333Y}) and (\ref{eq333F}) are equivalent if $\Phi$
is absolutely continuous, see e.g. Theorem IV.7.4 in \cite{Sa}, p.
119.\bigskip

\section{Connection with one more condition}

In this section we establish the useful connection of the conditions
of the Zorich--Lehto--Miklyukov-Suvorov type (\ref{eq3.333A})
further with one of the integral conditions from the last section.
\medskip

Recall that a function  $\Phi :[0,\infty ]\to [0,\infty ]$ is called
{\bf convex} if
$$
\Phi (\lambda t_1 + (1-\lambda) t_2)\ \le\ \lambda\ \Phi (t_1)\ +\
(1-\lambda)\ \Phi (t_2)$$ for all $t_1$ and $t_2\in[0,\infty ] $ and
$\lambda\in [0,1]$.\medskip

In what follows, ${\Bbb B}^n$ denotes the unit ball in the space
${\Bbb R}^n$, $n\ge 2$, \begin{equation}\label{eq5.5Cf} {\Bbb B}^n\
=\ \{\ x\in{\Bbb R}^n:\ |x|\ <\ 1\ \}\ .\end{equation}

\begin{lemma} \label{lem5.5C} Let $Q:{\Bbb B}^n\to [0,\infty ]$ be a measurable
function and let $\Phi:[0,\infty ]\to [0,\infty ]$ be a
non-decreasing convex function. Then \begin{equation}\label{eq3.222}
\int\limits_{0}^{1}\ \frac{dr}{rq(r)}\ \ge\ \frac{1}{n}\
\int\limits_{\lambda_nM}^{\infty}\ \frac{d\tau}{\tau \Phi^{-1}(\tau
)}
\end{equation} where $q(r)$ is the average of the function $Q(x)$
over the sphere $|x|=r$, \begin{equation}\label{eq555.555} M\ =\
\int\limits_{{\Bbb B}^n} \Phi (Q(x))\ dx\ ,\end{equation}
$\lambda_n=e/\Omega_n$ and $\Omega_n$ is the volume of the unit ball
in ${\Bbb R}^n$.\end{lemma}

\begin{remark}\label{rmk4.7vvv} In other words, \begin{equation}\label{eq3.555}
\int\limits_{0}^{1}\ \frac{dr}{rq(r)}\ \ge\ \frac{1}{n}\
\int\limits_{eM_n}^{\infty}\ \frac{d\tau}{\tau \Phi^{-1}(\tau )}
\end{equation} where \begin{equation}\label{eq333.333} M_n\colon =\
\frac{M}{\Omega_n}\ =\dashint\limits_{{\Bbb B}^n} \Phi (Q(x))\
dx\end{equation} is the mean value of the function $\Phi\circ Q$
over the unit ball. Recall also that by the Jacobi formula
$$
\Omega_n\ =\ \frac{\omega_{n-1}}{n}\ =\ \frac{2}{n}\cdot
 \frac{\pi^{\frac{n}{2}}}{\Gamma(\frac{n}{2})}\ =\
 \frac{\pi^{\frac{n}{2}}}{\Gamma(\frac{n}{2}+1)}\
$$ where $\omega_{n-1}$ is the area of the unit sphere in ${\Bbb R}^n$,
$\Gamma$ is the well--known gamma function of Euler,
$\Gamma(t+1)=t\Gamma(t)$. For $n=2$ we have that $\Omega_n=\pi$,
$\omega_{n-1}=2\pi$, and, thus, $\lambda_2=e/\Omega_2<1$.
Consequently, we have in the case $n=2$ that
\begin{equation}\label{eq777.777} \int\limits_{0}^{1}\
\frac{dr}{rq(r)}\ \ge\ \frac{1}{2}\ \int\limits_{M}^{\infty}\
\frac{d\tau}{\tau \Phi^{-1}(\tau )}\ .
\end{equation} In the general case we have that
$$
\Omega_{2m}\ =\ \frac{\pi^m}{m!}, \ \ \ \Omega_{2m+1}\ =\
\frac{2(2\pi)^m}{(2m+1)!!}\ ,
$$
i.e., $\Omega_n\to 0$ and, correspondingly, $\lambda_n\to\infty$ as
$n\to\infty$.\end{remark}

{\it Proof.} Note that the result is obvious if $M=\infty$. Hence we
assume further that $M<\infty$. Consequently, we may also assume
that $\Phi(t)<\infty$ for all $t\in [0,\infty)$ because in the
contrary case $Q\in L^{\infty}({\Bbb B}^n)$ and then the left hand
side in (\ref{eq3.222}) is equal to $\infty$. Moreover, we may
assume that $\Phi(t)$ is not constant (because in the contrary case
$\Phi^{-1}(\tau)\equiv\infty$ for all $\tau>\tau_0$ and hence the
right hand side in (\ref{eq3.222}) is equal to 0), $\Phi(t)$ is
(strictly) increasing, convex and continuous in a segment
$[t_*,\infty]$ for some $t_*\in [0,\infty)$ and
\begin{equation}\label{eq5.555Y} \Phi(t)\ \equiv\ \tau_0\ =\ \Phi(0)\
\ \ \ \ \ \ \ \forall\ t\ \in[0,t_*]\ .
\end{equation}

Next, setting \begin{equation}\label{eq5.555F}H(t)\ \colon =\ \log\
\Phi(t)\ ,\end{equation} we see by Proposition \ref{pr5.AAA} and
Remark \ref{rmk4.A} that
\begin{equation}\label{eq5.555G}H^{-1}(\eta)\ =\ \Phi^{-1}(e^{\eta})\
,\ \ \ \Phi^{-1}(\tau)\ =\ H^{-1}(\log\ \tau)\ .\end{equation} Thus,
we obtain that
\begin{equation}\label{eq5.555I}q(r) = H^{-1}\left(\log
\frac{h(r)}{r^n}\right) = H^{-1}\left(n\log \frac{1}{r} + \log\
h(r)\right)\ \ \ \ \ \ \ \ \forall\ r\ \in R_*\end{equation} where
$h(r)\ \colon =\ r^n\Phi(q(r))$ and $R_*\ =\ \{ r\in(0,1):\ q(r)\ >\
t_*\}$. Then also
\begin{equation}\label{eq5.555K}q(e^{-s})\ =\ H^{-1}\left(ns\ +\ \log\
h(e^{-s})\right)\ \ \ \ \ \ \ \ \forall\ s\ \in S_*\end{equation}
where $ S_*\ =\ \{ s\in(0,\infty):\ q(e^{-s})\ >\ t_*\}$.

Now, by the Jensen inequality \begin{equation}\label{eq5.555L}
\int\limits_0^{\infty} h(e^{-s})\ ds\ =\ \int\limits_0^{1} h(r)\
\frac{dr}{r}\ =\ \int\limits_0^{1} \Phi(q(r))\
r^{n-1}{dr}\end{equation}
$$\le\ \int\limits_0^{1}\left(\dashint_{S(r)}
\Phi(Q(x))\ d{\cal {A}}\right)\ r^{n-1}{dr}\ =\
\frac{M}{\omega_{n-1}}$$ where we use the mean value of the function
$\Phi\circ Q$ over the sphere $S(r)=\{ x\in{\Bbb R}^n: |x|=r\}$ with
respect to the area measure. Then
\begin{equation}\label{eq5.555N} |T|\ =\ \int\limits_{T}ds\ \le\
\frac{\Omega_n}{\omega_{n-1}}\ =\ \frac{1}{n}\end{equation} where
$T\ =\ \{\ s\in (0,\infty):\ \ \ h(e^{-s})\ >\ M_n\}$,
$M_n=M/{\Omega_n}$. Let us show that
\begin{equation}\label{eq5.555O}q(e^{-s})\ \le\ H^{-1}\left(ns\ +\ \log\
{M_n}\right)\ \ \ \ \ \ \ \ \ \ \forall\ s\in(0,\infty)\setminus
T_*\end{equation} where $T_*\ =\ T\cap S_*$. Note that
$(0,\infty)\setminus T_*  = [(0,\infty)\setminus S_*] \cup
[(0,\infty)\setminus T] = [(0,\infty)\setminus S_*] \cup
[S_*\setminus T]$. The inequality (\ref{eq5.555O}) holds for $s\in
S_*\setminus T$ by (\ref{eq5.555K}) because $H^{-1}$ is a
non-decreasing function. Note also that by (\ref{eq5.555Y})
\begin{equation}\label{eq5.K555}  e^{ns}{M_n} =
e^{ns} \dashint_{{\Bbb B}^n}\ \Phi(Q(x)) \ dx > \Phi(0) = \tau_0 \ \
\ \ \ \ \ \forall\ s\in(0,\infty)\ .\end{equation} Hence, since the
function $\Phi^{-1}$ is non-decreasing and
$\Phi^{-1}(\tau_0+0)=t_*$, we have by (\ref{eq5.555G}) that
\begin{equation}\label{eq5.L555} t_* < \Phi^{-1}\left({M_n}\
e^{ns}\right) = H^{-1}\left(ns\ +\ \log\ {M_n}\right)\ \ \ \ \ \ \ \
\forall\ s\in(0,\infty)\ .\end{equation} Consequently,
(\ref{eq5.555O}) holds for $s\in(0,\infty)\setminus S_*$, too. Thus,
(\ref{eq5.555O}) is true.\medskip

Since $H^{-1}$ is non--decreasing, we have by (\ref{eq5.555N}) and
(\ref{eq5.555O}) that \begin{equation}\label{eq5.555P}
\int\limits_{0}^{1}\ \frac{dr}{rq(r)}\ =\ \int\limits_{0}^{\infty}\
\frac{ds}{q(e^{-s})}\ \ge\ \int\limits_{(0,\infty)\setminus T_*}\
\frac{ds}{H^{-1}(ns + \Delta)}\ \ge\ \end{equation}
$$
\ge\ \int\limits_{|T_*|}^{\infty}\ \frac{ds}{H^{-1}(ns + \Delta)}\
\ge\ \int\limits_{\frac{1}{n}}^{\infty}\ \frac{ds}{H^{-1}(ns +
\Delta)}\ =\ \frac{1}{n}\int\limits_{1+{\Delta}}^{\infty}\
\frac{d\eta}{H^{-1}(\eta)}
$$
where $ \Delta=\log {M_n}$. Note that $1+\Delta = \log\ e{M_n}\ =
\log\ \lambda_nM$. Thus,
\begin{equation}\label{eq5.555S} \int\limits_{0}^{1}\
\frac{dr}{rq(r)}\ \ge\ \frac{1}{n}\int\limits_{\log \lambda_n{M}
}^{\infty}\ \frac{d\eta}{H^{-1}(\eta)}
\end{equation} and, after the replacement $\eta = \log\ \tau$, we
obtain (\ref{eq3.222}).

\begin{corollary} \label{cor3.1}
Let $Q:{\Bbb B}^n\to [0,\infty ]$ be a measurable function and let
$\Phi:[0,\infty ]\to [0,\infty ]$ be a non-decreasing convex
function. Then \begin{equation}\label{eq3.1} \int\limits_{0}^{1}\
\frac{dr}{rq^{\lambda}(r)}\ \ge\ \frac{1}{n}\
\int\limits_{\lambda_nM_*}^{\infty}\ \frac{d\tau}{\tau
\Phi^{-1}(\tau )}\ \ \ \ \ \ \ \ \ \ \forall\ \lambda\ \in\ (0,1)
\end{equation} where $q(r)$ is the average of the function $Q(x)$
over the sphere $|x|=r$, \begin{equation}\label{eq3.2} M_*\ =\
\int\limits_{{\Bbb B}^n} \Phi (Q_*(x))\ dx\ ,\end{equation} $Q_*$ is
the lower cut--off function of $Q$, i.e., $Q_*(x)=1$ if $Q(x)<1$ and
$Q_*(x)=Q(x)$ if $Q(x)\ge 1$.
\end{corollary}

Indeed, let $q_*(r)$ be the average of the function $Q_*(x)$ over
the sphere $|x|=r$. Then $q(r)\le q_*(r)$ and, moreover, $q_*(r)\ge
1$ for all $r\in (0,1)$. Thus, $q^{\lambda}(r)\le
q_*^{\lambda}(r)\le q_*(r)$  for all $\lambda\in (0,1)$ and hence by
Lemma \ref{lem5.5C} applied to the function $Q_*(x)$ we obtain
(\ref{eq3.1}).

\begin{theorem} \label{th5.555} Let $Q:{\Bbb B}^n\to [0,\infty ]$ be a measurable
function such that \begin{equation}\label{eq5.555}
\int\limits_{{\Bbb B}^n} \Phi (Q(x))\ dx\  <\ \infty\end{equation}
where $\Phi:[0,\infty ]\to [0,\infty ]$ is a non-decreasing convex
function such that
\begin{equation}\label{eq3.333a} \int\limits_{\delta_0}^{\infty}\ \frac{d\tau}{\tau
\Phi^{-1}(\tau )}\ =\ \infty \end{equation} for some $\delta_0\ >\
\tau_0\ \colon =\ \Phi(0).$ Then \begin{equation}\label{eq3.333A}
\int\limits_{0}^{1}\ \frac{dr}{rq(r)}\ =\ \infty \end{equation}
where $q(r)$ is the average of the function $Q(x)$ over the sphere
$|x|=r$.
\end{theorem}

\begin{remark}\label{rmk4.7www}  Note that (\ref{eq3.333a}) implies that
\begin{equation}\label{eq3.a333} \int\limits_{\delta}^{\infty}\ \frac{d\tau}{\tau
\Phi^{-1}(\tau )}\ =\ \infty \end{equation} for every $\delta\ \in\
[0,\infty)$ but (\ref{eq3.a333}) for some $\delta\in[0,\infty)$,
generally speaking, does not imply (\ref{eq3.333a}). Indeed, for
$\delta\in [0,\delta_0),$ (\ref{eq3.333a}) evidently implies
(\ref{eq3.a333})  and, for $\delta\in(\delta_0,\infty)$, we have
that
\begin{equation}\label{eq3.e333} 0\ \le\ \int\limits_{\delta_0}^{\delta}\
\frac{d\tau}{\tau \Phi^{-1}(\tau )}\ \le\
\frac{1}{\Phi^{-1}(\delta_0)}\ \log\ \frac{\delta}{\delta_0}\ <\
\infty
\end{equation} because $\Phi^{-1}$ is non-decreasing and
$\Phi^{-1}(\delta_0)>0$. Moreover, by the de\-fi\-ni\-tion of the
inverse function $\Phi^{-1}(\tau)\equiv 0$ for all $\tau \in
[0,\tau_0],$ $\tau_0=\Phi(0)$, and hence (\ref{eq3.a333}) for
$\delta\in[0,\tau_0),$ generally speaking, does not imply
(\ref{eq3.333a}). If $\tau_0 > 0$, then
\begin{equation}\label{eq3.c333} \int\limits_{\delta}^{\tau_0}\
\frac{d\tau}{\tau \Phi^{-1}(\tau )}\ =\ \infty\ \ \ \ \ \ \ \ \ \
\forall\ \delta\ \in\ [0,\tau_0)  \end{equation} However,
(\ref{eq3.c333}) gives no information on the function $Q(x)$ itself
and, consequently, (\ref{eq3.a333}) for $\delta < \Phi(0)$ cannot
imply (\ref{eq3.333A}) at all. \end{remark}

By  (\ref{eq3.a333}) the proof of Theorem \ref{th5.555} is reduced
to Lemma \ref{lem5.5C}.

\begin{corollary} \label{cor555} If $\Phi:[0,\infty ]\to [0,\infty ]$ is a
non-decreasing convex func\-tion and $Q$ satisfies the condition
(\ref{eq5.555}), then every of the conditions
(\ref{eq333Y})--(\ref{eq333A}) implies (\ref{eq3.333A}).\medskip

Moreover, if in addition $\Phi(1)<\infty$ or $q(r)\ge 1$ on a subset
of $(0,1)$ of a positive measure, then
\begin{equation}\label{eq3.3} \int\limits_{0}^{1}\
\frac{dr}{rq^{\lambda}(r)}\ =\ \infty\ \ \ \ \ \ \ \ \ \forall\
\lambda\ \in\ (0,1)
\end{equation}
and also
\begin{equation}\label{eq3.3AB} \int\limits_{0}^{1}\
\frac{dr}{r^{\alpha}q^{\beta}(r)}\ =\ \infty\ \ \ \ \ \ \ \ \
\forall\ \alpha\ge 1 ,\ \beta\ \in\ (0,\alpha]
\end{equation}
\end{corollary}

{\it Proof}. First of all, by Theorems \ref{pr4.1aB} and
\ref{th5.555} every of the conditions (\ref{eq333Y})--(\ref{eq333A})
implies (\ref{eq3.333A}). Now, if $q(r)\ge 1$ on a subset of $(0,1)$
of a positive measure, then also $Q(x)\ge 1$ on a subset of ${\Bbb
B}^n $ of a positive measure and, consequently, in view of
(\ref{eq5.555}) we have that $\Phi(1)<\infty$. Then
\begin{equation}\label{eq555} \int\limits_{{\Bbb B}^n} \Phi
(Q_*(x))\ dx\  <\ \infty\end{equation} where $Q_*$ is the lower cut
off function of $Q$ from  Corollary \ref{cor3.1}. Thus, by Corollary
\ref{cor3.1} and Remark \ref{rmk4.7www} we obtain that every of the
conditions (\ref{eq333Y})--(\ref{eq333A}) implies (\ref{eq3.3}).
Finally, (\ref{eq3.3AB}) follows from (\ref{eq3.3}) by the Jensen
inequality.

\begin{remark}\label{rmkvvv}
Note that if we have instead of (\ref{eq5.555}) the condition
\begin{equation}\label{eq5.555A}
\int\limits_{D} \Phi (Q(x))\ dx\  <\ \infty\end{equation} for some
measurable function $Q: D\to [0,\infty]$ given in a domain $D\subset
{{\Bbb R}^n}$, $n\ge 2$, then also \begin{equation}\label{eq5.555B}
\int\limits_{|x-x_0|<r_0} \Phi (Q(x))\ dx\  <\ \infty\end{equation}
for every $x_0\in D$ and $r_0<\mbox{dist}\, (x_0, \partial D)$ and
by Theorem \ref{th5.555} and Corollary \ref{cor555}, after the
corresponding linear replacements of variables, we obtain that
\begin{equation}\label{eq3.3A} \int\limits_{0}^{r_0}\
\frac{dr}{rq_{x_0}^{\lambda}(r)}\ =\ \infty\ \ \ \ \ \ \ \ \
\forall\ \lambda\ \in\ (0,1]
\end{equation} where $q_{x_0}(r)$ is the average of the function $Q(x)$
over the sphere $|x-x_0|=r$.

If $D$ is a domain in the extended space $\overline{{\Bbb
R}^n}={{\Bbb R}^n}\cup\{\infty\}$ and $\infty\in D$, then in the
neighborhood $|x|>R_0$ of $\infty$ we may use the condition
\begin{equation}\label{eq5.555C} \int\limits_{|x|>R_0} \Phi
(Q(x))\ \frac{dx}{\ \ |x|^{2n}}\  <\ \infty\end{equation} that is
equivalent to the condition \begin{equation}\label{eq5.555D}
\int\limits_{|x|<r_0} \Phi (Q'(x))\ dx\  <\ \infty\end{equation}
where $r_0=1/R_0$ and $Q'(x)=Q(x/|x|^2)$, i.e. $Q'(x)$ is obtained
from $Q(x)$ by the inversion of the independent variable $x\to
x/|x|^2$, $\infty \to 0$, with respect to the unit sphere $|x|=1$.

Thus, by Theorem \ref{th5.555} and Corollary \ref{cor555} the
condition (\ref{eq5.555C}) imply the equality
\begin{equation}\label{eq3.3B} \int\limits_{R_0}^{\infty}\
\frac{dR}{Rq_{\infty}^{\lambda}(R)}\ =\ \infty\ \ \ \ \ \ \ \ \
\forall\ \lambda\ \in\ (0,1]
\end{equation} where $q_{\infty}(R)$ is the average of $Q$
over the sphere $|x|=R$.

Finally, if $D$ is an unbounded domain in ${{\Bbb R}^n}$ or a domain
in $\overline{{\Bbb R}^n}$ (\ref{eq5.555A}) should be replaced by
the following condition \begin{equation}\label{eq5.555E}
\int\limits_{D} \Phi (Q(x))\ dS(x)\  <\ \infty\end{equation} where
$dS(x)= dx/(1+|x|^2)^n$ is of a cell of the spherical volume. Here
the spherical distance \begin{equation} \label{eq1.5a}
 s(x,y)\  =\
\frac{|x-y|}{(1+|x|^2)^{\frac{1}{2}}(1+|y|^2)^{\frac{1}{2}}}\ \ \
{\mbox{if}}\ \ \  x\ \ne\ \infty\ne y\ ,\end{equation}
$$ s(x,\infty )\ =\ \frac{1}{(1+|x|^2)^{\frac{1}{2}}}\ \ \
{\mbox{if}}\ \ \  x\ \ne\ \infty\ .$$ It is easy to see that
$$dS(x)\ge (1+\rho^2)^{-n}dx$$ in every bounded part of $D$ where
$|x|<\rho$ and $$dS(x)\ge 2^{-n}\frac{dx}{\ \ |x|^{2n}}$$ in a
neighborhood of $\infty$ where $|x|\ge 1$. Hence the condition
(\ref{eq5.555E}) implies (\ref{eq5.555B}) as well as
(\ref{eq5.555C}). Thus, under at least one of the conditions
(\ref{eq333F})--(\ref{eq333A}) the condition (\ref{eq5.555E})
implies (\ref{eq3.3A}) and (\ref{eq3.3B}) if the function $\Phi$ is
convex and non-decreasing.

\bigskip

Recently it was established that the conditions
(\ref{eq333F})--(\ref{eq333A}) are not only sufficient but also
necessary for the degenerate Beltrami equations with integral
constraints of the type (\ref{eq5.555}) on their
cha\-rac\-te\-ris\-tics to have ho\-meo\-mor\-phic solutions of the
class $W^{1,1}_{loc}$, see \cite{BJ$_2$}, \cite{GMSV}, \cite{IM$_2$}
and \cite{RSY$_5$}.\bigskip

Moreover, the above results will have significant corollaries to the
local and boundary behavior of space mappings in various modern
classes with integral constrains for dilatations, see e.g.
\cite{MRSY}. \end{remark}

\medskip
\noindent {\bf Vladimir Ryazanov:}\\ Institute of Applied
Mathematics\\ and Mechanics, NAS of Ukraine, \\ ul. Roze Luxemburg
74,\\ 83114, Donetsk, UKRAINE\\ Email: {\tt
vlryazanov1@rambler.ru}\\
\medskip

\noindent {\bf Uri Srebro:}\\ Technion - Israel Institute of
Technology, \\ Haifa 32000, ISRAEL\\ Fax: 972--4--8324654\\
Email: {\tt srebro@math.technion.ac.il}\\

\medskip
\noindent {\bf Eduard Yakubov:}\\ Holon Institute of Technology,\\
52 Golomb St., P.O.Box 305,\\ Holon 58102,
ISRAEL\\ Email: {\tt yakubov@hit.ac.il}\\

\end{document}